\documentclass{article} 
\usepackage{nips15submit_e,times}
\usepackage{url}
\usepackage{graphicx}
\usepackage{caption}
\usepackage{subcaption}

\usepackage{epsfig}
\usepackage{epstopdf} 

\usepackage{comment}
\usepackage{amsmath,amssymb, amsthm}

\usepackage{fullpage}

\usepackage{multirow}
\usepackage{amsfonts}
\usepackage{layout}
\usepackage{url}
\usepackage{color}
\usepackage{graphicx}
\usepackage{algorithmic}
\usepackage{algorithm}
\usepackage{verbatim}
\usepackage{epsfig}
\usepackage{xcolor}
\usepackage{tikz}

\PassOptionsToPackage{normalem}{ulem}
\usepackage{ulem}

\usepackage{cases}

\providecolor{added}{rgb}{0,0,1}
\providecolor{deleted}{rgb}{1,0,0}

\definecolor{orange}{RGB}{255,127,0}

\newcommand{\R}{\mathbb{R}}

\newcommand{\E}{\mathbf{E}}

\newcommand{\vc}[2]{#1^{(#2)}}

\newcommand{\m}[1]{~\mbox{#1}~}

\newcommand{\norm}[1]{\|#1\|}

\DeclareMathOperator{\Exp}{\mathbf{E}}           

\newtheorem{lemma}{Lemma}

\newtheorem{theorem}{Theorem}

\theoremstyle{plain}

\theoremstyle{definition}

\newtheorem{definition}{Definition}

\newcommand*{\starnr}{\stepcounter{equation}\tag{\theequation}}

\usepackage{hyperref}

\title{Dual Free SDCA for Empirical Risk Minimization with Adaptive Probabilities}

\author{
    Xi He \\
    Industrial and Systems Engineering\\
    Lehigh University, USA\\
    \texttt{xih314@lehigh.edu} \\
    \And
    Martin Tak\'a\v{c} \\
    Industrial and Systems Engineering\\
    Lehigh University, USA \\
    \texttt{takac.mt@gmail.com} \\
}

\nipsfinalcopy 

\begin{document}

\maketitle

\begin{abstract}
    In this paper we develop dual free SDCA
    with adaptive probabilities for regularized empirical risk minimization. This extends recent work of Shai Shalev-Shwartz [SDCA without Duality, arXiv:1502.06177]
    to allow non-uniform selection of "dual" coordinate in SDCA. Moreover, the probability can change over time, making it more efficient than uniform selection. Our work focuses on generating adaptive probabilities through iterative process, preferring to choose coordinate with highest potential to decrease sub-optimality.  We also propose a practical variant Algorithm adfSDCA\texttt{+} which is more aggressive. The work is concluded with multiple experiments which  shows   efficiency of proposed algorithms.
\end{abstract}

\section{Introduction}
We study the $\ell_2$-regularized Empirical Risk Minimization (ERM) problem, which is widely used in the field of  machine learning. Given training examples $(x_1,y_1), \dots,  (x_n,y_n) \in \R^d \times \R$, loss functions $\phi_1, \dots, \phi_n: \R \rightarrow \R$ and a regularization parameter $\lambda > 0$, $\ell_2$-regularized ERM problem is an optimization problem of the form
\begin{equation}\label{Prob: L2EMR}
\tag{P}
\min_{w\in \R^d} P(w) := \tfrac{1}{n}\textstyle{\sum}_{i=1}^n\phi_i(w^Tx_i) + \tfrac{\lambda}{2}\norm{w}^2,
\end{equation}
where the first part of the objective function is the \emph{data fitting term}
and the second part it the regularization which prevents over-fitting.

Various methods were proposed for solving this problem over past few years,
many of them are trying to handle the problem \eqref{Prob: L2EMR} directly
including \cite{shalev2011pegasos,johnson2013accelerating,schmidt2013minimizing,defazio2014saga,
roux2012stochastic,nitanda2014stochastic,konevcny2014ms2gd}, others are trying to solve its dual formulation \cite{hsieh2008dual}:
\begin{equation}
\tag{D} \label{Prob:dual}
\max_{\alpha \in \R^n} D(\alpha) :=
-\frac1n \textstyle{\sum}_{i=1}^n \phi^*_i(-\alpha_i) -\tfrac{\lambda}{2} \|\tfrac1{\lambda n} X^T \alpha\|^2,
\end{equation}
where $X^T=[x_1,\dots, x_n]\in \R^{d \times n}$ is the data matrix and
$\phi_i^*$ is a convex conjugate function of $\phi_i$.

One of the most popular method for solving \eqref{Prob:dual} is so-called Stochastic Dual Coordinate Ascent (SDCA).
In each iteration $t$ of SDCA,
a coordinate $i\in\{1, \dots, n\}$ is chosen uniformly at random and current iteration $\vc{\alpha}{t}$ is updated to $\vc{\alpha}{t+1}:=\vc{\alpha}{t}+  \delta^* e_i$, where
$\delta^* = \arg\max_{\delta\in \R} D(\vc{\alpha}{t}+ \delta e_i)$.
There has been a lot of work done for analysing the complexity theory of SDCA under various assumptions imposed on $\phi_i^*$ including a pioneering work of Nesterov \cite{nesterov2012efficiency} and others \cite{richtarik2014iteration,
tappenden2015complexity,
takavc2015distributed,
takavc2013mini}.

One fascinating algorithmic change which turned out to improve the convergence in practice is so-called \emph{importance sampling}, i.e. a step when we sample coordinate $i$ with arbitrary probability $p_i$ \cite{zhao2014stochastic,
csiba2015stochastic,
csiba2015primal,schmidt2015non,
richtarik2013optimal}.
This turns out to outperform the na\"ive uniform selection and in some cases help to decrease the number of iterations needed to achieve a desired accuracy by few folds.

\paragraph{Assumptions.}
In this work we assume that $\forall i\in \{1,\dots,n\}:=[n]$, the loss function $\phi_i$ is $\tilde L_i$-smooth with $\tilde L_i>0$, i.e.  for any given $\beta, \delta \in \R$, we have
\begin{equation}
    | \phi_i'(\beta) - \phi_i'(\beta+\delta)| \leq \tilde L_i|\delta |.
\end{equation}
It is a simple observation that
also function $\phi_i(x_i^T \cdot  ): \R^d \to \R$ is $L_i$ smooth, i.e. $\forall i \in [n], \forall w,\tilde w\in \R^d$
there exists a constant  $L_i\leq \|x_i\|^2 \tilde L_i $ such that
\begin{equation}
    \|\nabla\phi_i( x_i^T w)-\nabla\phi_i(x_i^T\bar{w})\| \leq L_i\norm{w-\bar{w}}, \m{for all}.
\end{equation}
We will denote by  $L=\max_i L_i$.




\subsection{Contributions} 
\label{ssec:contributions}
In this work we modify dual free SDCA proposed by Shalev-Shwartz \cite{DBLP:journals/corr/Shalev-Shwartz15} to allow adaptive adjustment of probabilities in non-uniform selection of coordinate. Note that the method is dual free, and hence in contrast of classical SDCA, when the update is trying to maximize the dual objective \eqref{Prob:dual}, we define the update differently (see Section \ref{sec:algorithm_adfSDCA} for more details).

Allowing adaptive non-uniform selection of coordinate leads  to efficient utilization of computational resource and the algorithm achieves a better complexity bound than \cite{DBLP:journals/corr/Shalev-Shwartz15}.
We showed that the error after $T$ iterations is decreased by factor of
$\prod_{t=1}^T(1-\tilde{\theta}_t)$.
Here, $1-\tilde{\theta}_t \in (0,1)$ is a  parameter which depends on current iterate $\vc{\alpha}{t}$ and probability distribution over the choice of coordinate to be used in next iteration.
By changing the strategy from uniform selection to adaptive, we are making each $1-\tilde{\theta}_t$ smaller, hence improving the convergence rate.


\section{The Adaptive Dual Free SDCA Algorithm} 
\label{sec:algorithm_adfSDCA}

In dual free SDCA as proposed by \cite{DBLP:journals/corr/Shalev-Shwartz15}
we maintain two sequence of iterates: $\{\vc{w}{t}\}_{t=0}^\infty$ and $\{\vc{\alpha}{t}\}_{t=0}^\infty$.
The updates in the algorithm are done in such a way that the well known
primal-dual relation mapping holds for every $t$:
\begin{equation}\label{asdffsafsafa}
\vc{w}{t} = \frac{1}{\lambda n}\sum_{i=1}^n \vc{\alpha_i}{t} x_i.
\end{equation}
In Algorithm \ref{Alg: adfSDCA} (heuristic = False)
we start with some initial solution $\vc{\alpha}{0}$ and
we define $\vc{w}{0}$ using \eqref{asdffsafsafa}.
In each iteration of  Algorithm \ref{Alg: adfSDCA} we compute the dual residuals $\vc{\kappa}{t}$ and generate a probability distribution $\vc{p}{t}$ based on the residuals.
Afterwards, we pick a coordinate $i \in [n]$ using the generated probability distribution and we take the step by updating $i$-th coordinate of $\alpha$ and making a corresponding update to $w$ such that \eqref{asdffsafsafa}
will be preserved.
\begin{algorithm}[h!]
    \caption{Adaptive Dual Free SDCA (adfSDCA)}
    \label{Alg: adfSDCA}
    \begin{algorithmic}[1]
        \STATE {\bf Input:} Data: $\{x_i, \phi_i\}_{i=1}^n$
        \STATE  {\bf Initialization:} Choose  $ \alpha^{(0)} \in \R^n$
        \STATE Set $w^{(0)} = \frac{1}{\lambda n}\sum_{i=1}^n \alpha_i^{(0)}x_i $
        \FOR {$t=0,1,2,\dots$}
        \IF {\m{heuristic} $\&\&\mod(t, n) == 0$}
        \STATE Calculate dual residue $\kappa^{(t)}_i = \phi'_i(x_i^Tw^{(t)}) + \alpha^{(t)}_i$, for all $i\in [n]$
        \STATE Generating adapted probabilities distribution $p^{(t)} \sim \kappa^{(t)}$
        \ELSIF {!\m{heuristic}}
        \STATE Calculate dual residue $\kappa^{(t)}_i = \phi'_i(x_i^Tw^{(t)}) + \alpha^{(t)}_i$, for all $i\in [n]$
        \STATE Generating adapted probabilities distribution $p^{(t)} \sim \kappa^{(t)}$
        \ENDIF
        \STATE Select coordinate $i$ from $[n]$ according to $p^{(t)}$
        \STATE {\bf Update:} $\alpha_i^{(t+1)} = \alpha_i^{(t)} - \theta (p_i^{(t)})^{-1} \kappa^{(t)}_i$
        \STATE \textbf{Update: }$w^{(t+1)} = w^{(t)} - \theta(n \lambda p_i^{(t)})^{-1} \kappa^{(t)}_ix_i$

        \STATE {\bf if} \m{heuristic} {\bf then} \textbf{Update: } $p^{(t+1)}_i = p^{(t)}_i/s$
        \ENDFOR
    \end{algorithmic}
\end{algorithm}
Note that if for some $i$ we have $\kappa_i = 0$,
this means that  $\alpha_i=-\phi'_i(w^Tx_i)$, which indicates that
if we would choose coordinate $i$ in given iteration,
both $\alpha$ and $w$
would be unchanged.
On the other hand, large value of $|\kappa_i|$
indicates that the step which is going to be taken will be very large, hoping to improve sub-optimality of current solutions.

\begin{definition} (Coherence \cite{csiba2015stochastic})\label{Def: coherence}
Probability vector  $p\in \R^n$ is coherent with dual residue $\kappa\in \R^n$ if for any index $i$ in the support set of $\kappa$, denoted by $I_ \kappa := \{i \in [n]: \kappa_i \neq 0\}$, we have $p_i >0$ and $p_i = 0$ otherwise. We use $p\sim \kappa$ to represent this coherent relation.
\end{definition}




\paragraph{Adaptive dual free SDCA as reduced variance SGD method.}
Reduced variance SGD methods have became very popular over the last few years
\cite{konevcny2013semi,johnson2013accelerating,roux2012stochastic,defazio2014saga}.
It was show in  \cite{DBLP:journals/corr/Shalev-Shwartz15}
that uniform dual free SDCA is an instance of reduced variance SGD algorithm
(the variance of the stochastic gradient can be bounded by some measure of sub-optimality of current iterate).
Note that conditioned on $w^{(t-1)}$, we have
\begin{align}\label{eq: nabla(P)}
    \E[w^{(t)}] & = w^{(t-1)} - \E\left[\frac{\theta}{n \lambda p_i}(\phi'_i(x_i^Tw^{(t)})+\alpha_i^{(t)})^T x_i\right]    
    = w^{(t-1)} - \frac{\theta}{\lambda}\nabla P(w^{(t-1)}) .
\end{align}

Therefore, Algorithm \ref{Alg: adfSDCA} is eventually a variant of Stochastic Gradient Descent method. However, we can prove that the variance of the update goes to zero as we converge to an optimum, which is not true for vanilla Stochastic Gradient Descent.
Similarly as in  \cite{DBLP:journals/corr/Shalev-Shwartz15}
we can show that $\Exp[ (\frac{1}{p_i} \vc{\kappa_i}{t})^2 ] $
can be bound by sub-optimality of point $\vc{\alpha}{t}$.

\section{Convergence analysis} 
\label{sec:convergence_analysis}
In this section, we state the main convergence results.
We will limit ourself only to the case when each $\phi_i$ is convex, but this assumption can be relaxed and it is enough to assume that the average of $\phi_i$ is convex (however, the result will be a bit worse).

\begin{theorem}\label{thm: seperate convex}
Assume that for each $i \in [n]$ $\phi_i$ is $ L$-smooth and convex, then for
any iteration $t$ following holds
\begin{align}\label{eq: separate convex conclusion}
&\E[\frac{1}{n}\norm{\alpha^{(t+1)}-\alpha^*}^2 + \gamma  \norm{w^{(t+1)}-w^*}^2] - (1- \theta^{(t)})(\frac{1}{n}\norm{\alpha^{(t)}-\alpha^*}^2 + \gamma  \norm{w^{(t)}-w^*}^2)\\
&\leq \sum_{i=1}^n(-\frac{\theta^{(t)}}{n}(1- \frac{\theta^{(t)}}{p_i})+\frac{(\theta^{(t)})^2v_i \gamma }{n^2 \lambda^2 p_i})(\kappa_i^{(t)})^2,\nonumber
\end{align}
where $\gamma  = \lambda L$, $v_i = \norm{x_i}^2$ and $\kappa_i^{(t)}$ is residue of coordinate $i$ at $t$-th iteration.
\end{theorem}
Note that if the right hand side of   \eqref{eq: separate convex conclusion} is negative, we can obtain
\begin{equation}
    \E[\frac{1}{n}\norm{\alpha^{(t+1)}-\alpha^*}^2 + \gamma\norm{w^{(t+1)}-w^*}^2] \leq (1- \theta^{(t)})(\frac{1}{n}\norm{\alpha^{(t)}-\alpha^*}^2 + \gamma  \norm{w^{(t)}-w^*}^2).
\end{equation}
To guarantee their negativity, we can use any $\theta^{(t)}$ that is less than the function $\Theta(\cdot, \cdot): \R_+^n \times \R_+^n \rightarrow \R$ defined as
$
\Theta(\kappa, p) := \frac{n \lambda^2 \sum_{i\in I_ \kappa}\kappa_i^2}{\sum_{i\in I_ \kappa}(v_i \gamma+n \lambda^2)p_i^{-1}\kappa_i^2}.
$
The larger  the $\theta$ the better progress our algorithm will make. The optimal probability which will allow the highest $\theta$ can be obtained by solving following optimization problem:
\begin{align}\label{Prob: to derive best probabilities}
\max _{p\in \R^n_+: \sum_{i\in I_ \kappa} p_i=1} \quad&\Theta(\kappa, p).
\end{align}
It turns out that one can derive the optimal solution in a closed form and is given by following Lemma.
\begin{lemma}\label{lem: optimal probabilities}
    The optimal solution $p^*(\kappa)$ of \eqref{Prob: to derive best probabilities} is
    \begin{equation}\label{eq: optimal probabilities}
    p^*_i(\kappa) = \tfrac{\sqrt{v_i \gamma+n \lambda^2}|\kappa_i|}{\sum_{i\in I_ \kappa}\sqrt{v_i \gamma+n \lambda^2}|\kappa_i|}.
\end{equation}
\end{lemma}




Comparing with conclusion in \cite{csiba2015primal}, their results are weaker since they allow any fixed sampling distribution $p$. Out result enjoys better convergence rate at each iteration by setting sampling probabilities as \eqref{eq: optimal probabilities}. The shortage is at each iteration, we need extra $O(nnz(\{x_1,\dots,x_n\}))$ operations to derive $\kappa$.
Comparing with conclusion in \cite{csiba2015stochastic}, their optimal probabilities can be only applied on quadratic loss function. While in our case, it can apply on any convex loss function and specific non-convex function (average convex).

In order to overcome the shortage we made above, in this paper we show one other heuristic approach to apply adaptive probabilities. The motivation behind it is as follows: once we update one coordinate, it is natural that dual residual at this coordinate will decrease. Instead of calculating the  dual residuals at next coordinate to derive adaptive probabilities, we simply shrink the probability $p_i$ when $i$ was the last coordinate which was updated (see Algorithm \ref{Alg: adfSDCA} when heuristic = True). Obviously, this is not a exact algorithm, but still, we can get a relative much better practical results (see numerical experiments). We refer to this algorithm as adfSDCA+.
\section{Numerical experiments} 
\label{sec:numerical_experiments}

In this section we will compare the adfSDCA with its uniform variant dfSCDA \cite{DBLP:journals/corr/Shalev-Shwartz15} and also with Prox-SDCA \cite{shalev2014accelerated}. We used two loss functions, quadratic loss $ \phi_i(w^Tx_i) = \frac{1}{2}(w^Tx_i-y_i)$
and logistic loss $  \phi_i(w^Tx_i) = \log(1+\exp(-y_iw^Tx_i))$.
We did our experiments on standard datasets rcv1: ($n=20,242; d=47,237$), 
and mushrooms: $(n=8,124; d=112)$.

\begin{figure}[!ht]
    \centering
    \includegraphics[width=0.3\textwidth]{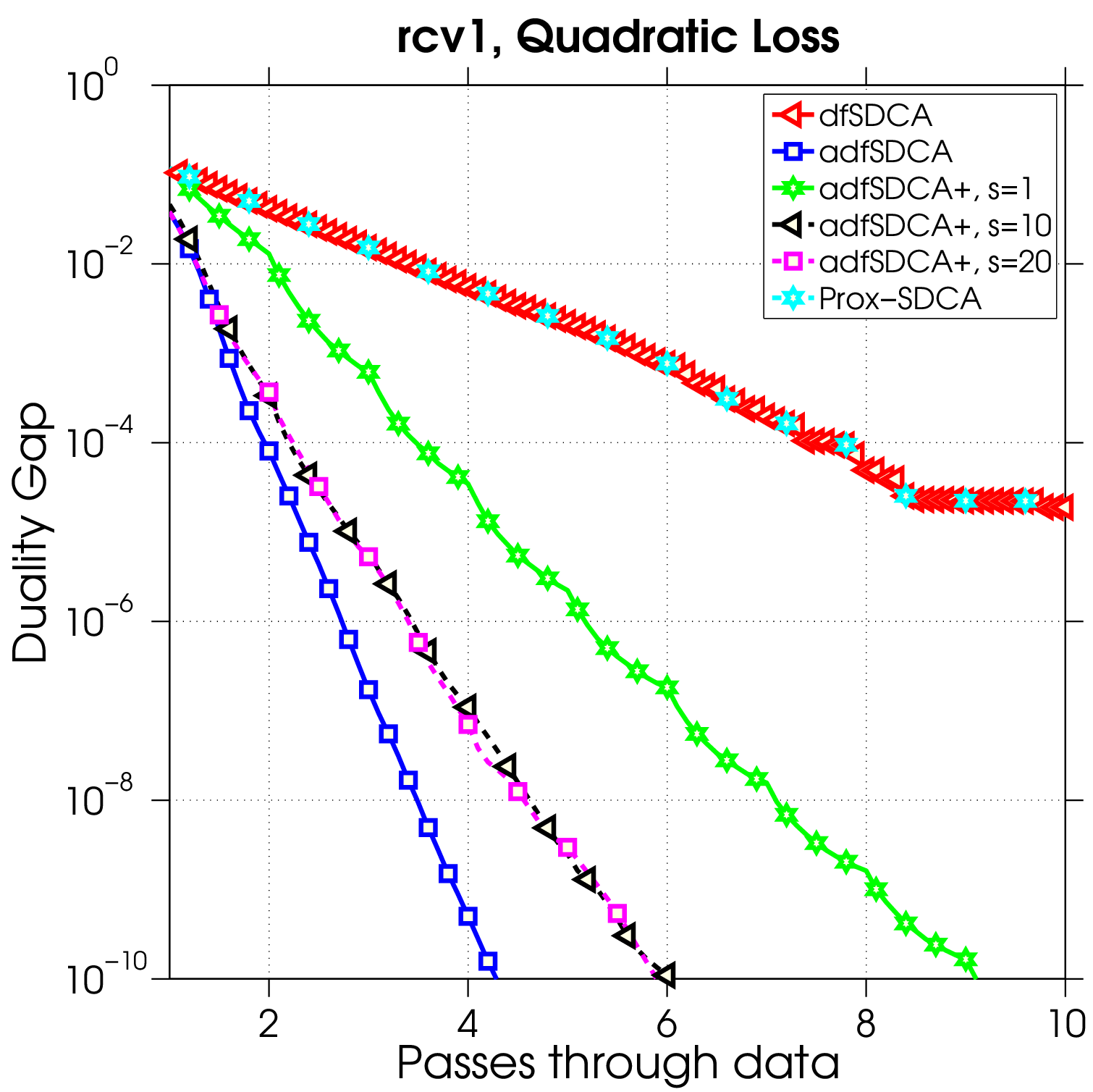}
    \includegraphics[width=0.3\textwidth]{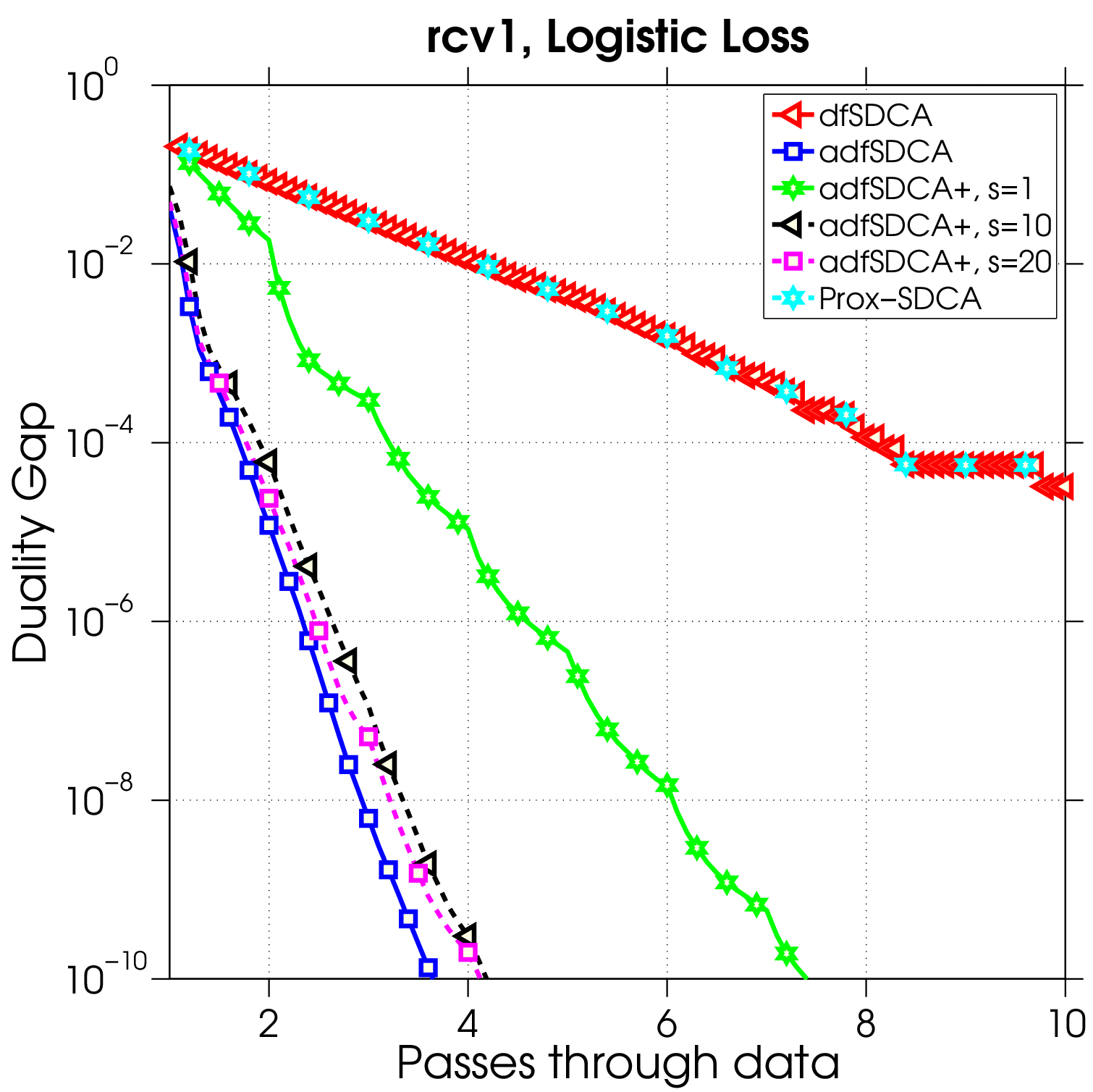}
    \includegraphics[width=0.3\textwidth]{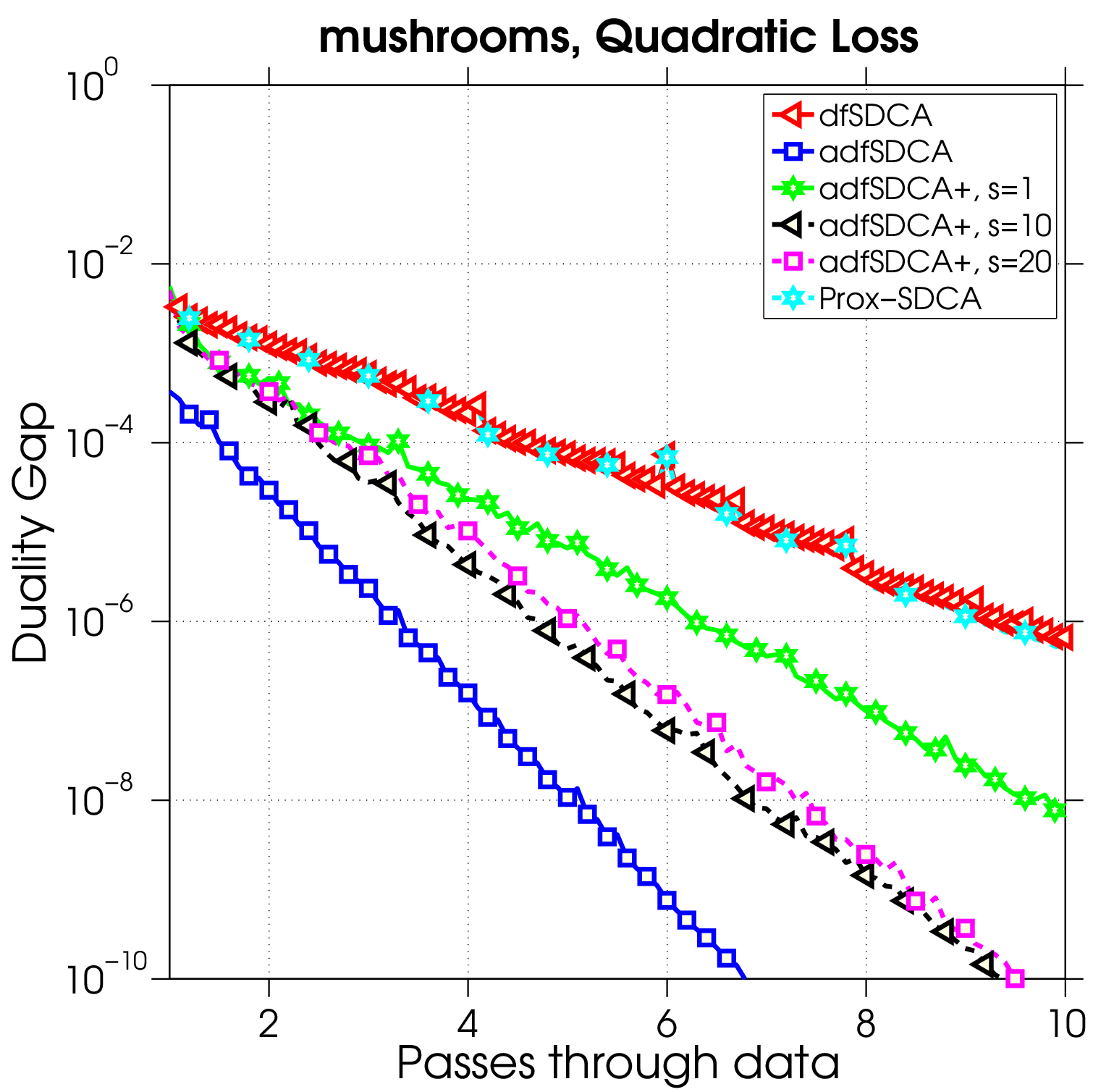}
    \caption{Comparing number of iterations among various algorithms.}
    \label{Fig:evolution}
\end{figure}

\begin{figure}[!ht]
    \centering
    \includegraphics[width=0.3\textwidth]{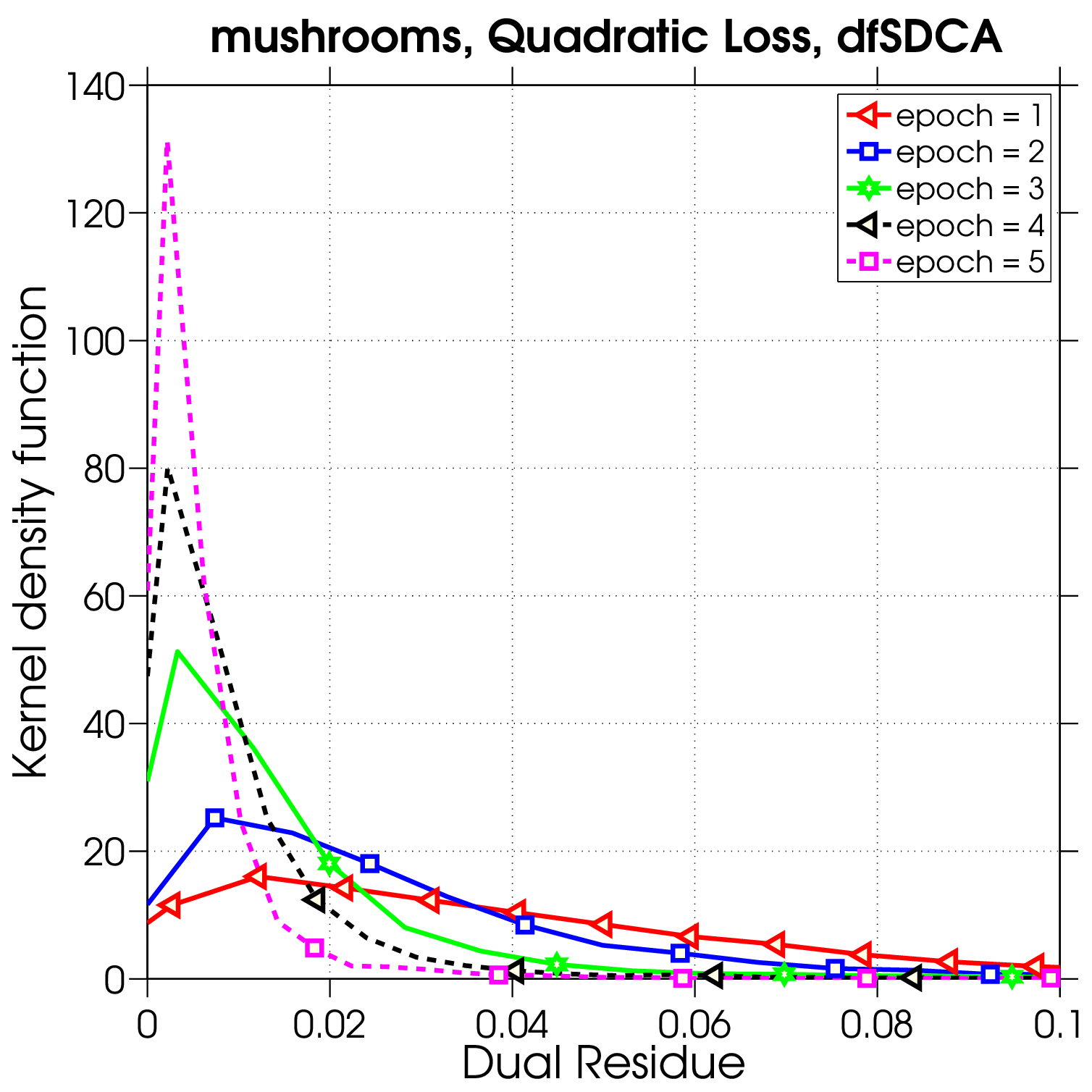}
    \includegraphics[width=0.3\textwidth]{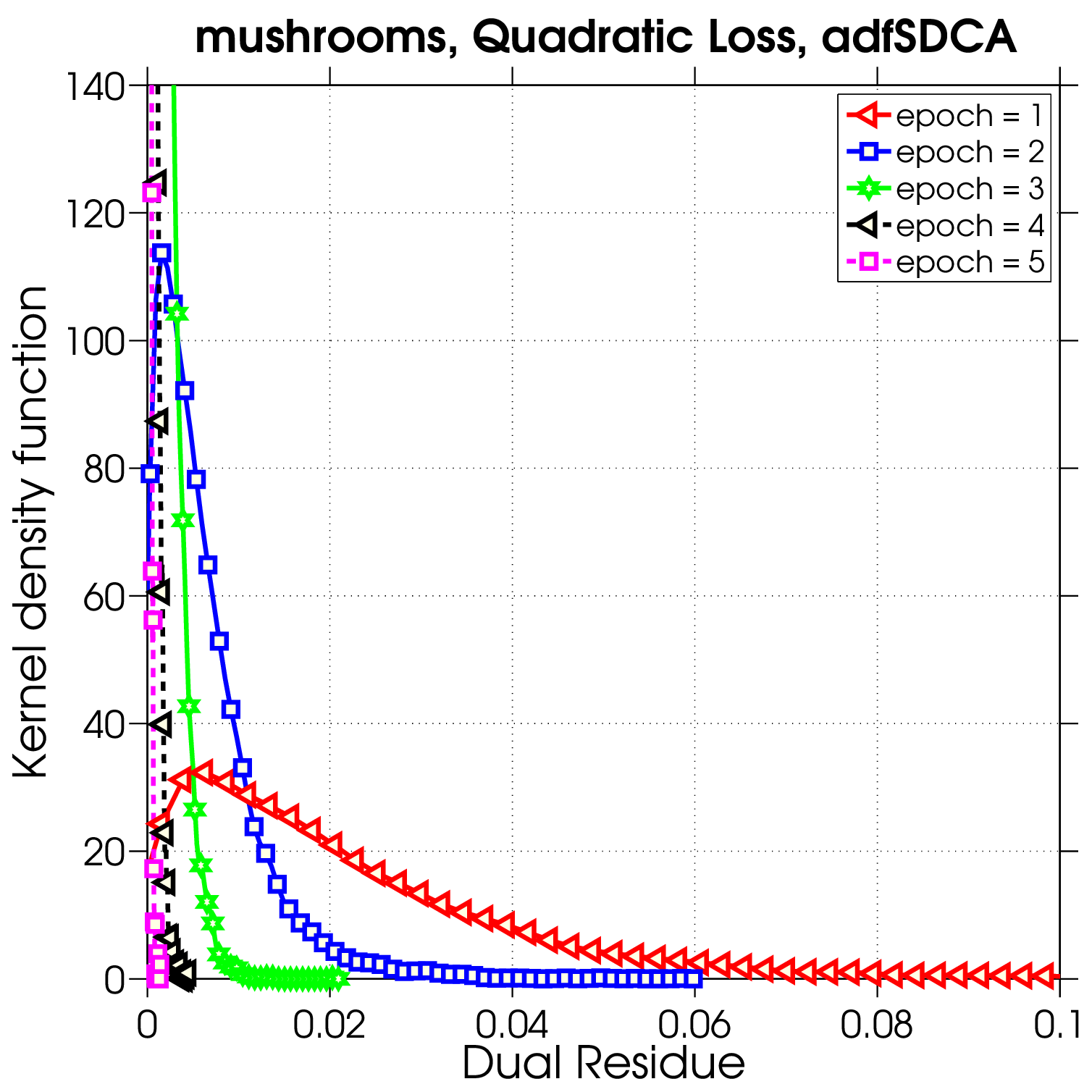}
    \caption{Comparing absolute value of dual residuals at each epoch between dfSDCA and adfSDCA.}
    \label{Fig:residuals}
\end{figure}

Figure  \ref{Fig:evolution} compares the evolution of duality gap for
various versions of our algorithm and shows the 2 state-of-the-art algorithms.
In this case all of our variants are out-performing the dfSDCA and Prox-SDCA algorithms. \\
Figure \ref{Fig:residuals} shows the estimated density function
of $|\vc{\kappa}{t}|$ after 1,2,3,4,5 epochs for uniform dfSDCA and our adaptive variant adfSDCA. As one can observe, the adaptive scheme is pushing the high residuals towards zero much faster. E.g. after 2 epochs, almost all residuals are below $0.03$ for adfSDCA case, whereas the uniform dfSDCA has still many residuals above $0.06$.

\bibliographystyle{plain}
\bibliography{ref}

\clearpage
\begin{comment}
\section*{Appendix} 
\label{sec:appendix}
To simplify proofs, we introduce
\begin{align}
    A^{(t)} = \frac{1}{n}\norm{\alpha^{(t)}-\alpha^*}^2 \m{and} B^{(t)} = \norm{w^{(t)}-w^*}^2,
\end{align}

where $w^*$  is an optimum of \eqref{Prob: L2EMR} and we define $\alpha_i^* = -\phi'_i(x_i^Tw^*)$. We can then derive that $w^* = \frac{1}{\lambda n}\sum_{i=1}^n\alpha_i^*x_i$ since  $\nabla P(w^*) = \frac{1}{n}\sum_{i=1}^n\phi'_i(x_i^Tw^*)x_i + \lambda w^* = 0$ at optimum $w^*$.

We claim the following two useful lemmas at first.
\begin{lemma}
    \begin{align}
        \E[A^{(t+1)} - A^{(t)}] &= -\theta A^{(t)} + \frac{\theta}{n}\sum_{i=1}^n((-\phi'_i(x_i^Tw^{(t)})-\alpha_i^*)^2 - (1- \frac{\theta}{p_i})(\kappa^{(t)}_i)^2),\label{eq: E[A]}\\
        \E[B^{(t+1)} - B^{(t)}] & = -\frac{2 \theta }{\lambda}\nabla P(w^{(t)})^T(w^{(t)}-w^*) + \sum_{i=1}^n \frac{\theta^2v_i(\kappa^{(t)}_i)^2}{n^2 \lambda^2p_i}, \m{where } v_i = \norm{x_i}^2.\label{eq: E[B]}
    \end{align}
\end{lemma}

\begin{proof}
    Note that
    \begin{align}
        &A^{(t+1)} - A^{(t)} \\
        =& \frac{1}{n}\norm{\alpha^{(t+1)}-\alpha^*}^2 - \frac{1}{n}\norm{\alpha^{(t)}-\alpha^*}^2\\
        =&\frac{1}{n}(\alpha_i^{(t)} - \frac{\theta}{p_i}(\phi'_i(x_i^Tw^{(t)})+\alpha_i^{(t)})-\alpha_i^*)^2 - \frac{1}{n}(\alpha_i^{(t)} - \alpha_i^*)^2\\
        =&\frac{1}{n}[(1- \frac{\theta}{p_i})(\alpha_i^{(t)}-\alpha_i^*)^2+\frac{\theta}{p_i}(-\phi'_i(x_i^Tw^{(t)})-\alpha_i^*)^2-(1- \frac{\theta}{p_i})\frac{\theta}{p_i}(\alpha_i^t+\phi'_i(x_i^Tw^{(t)}))^2] - \frac{1}{n}(\alpha_i^{(t)} - \alpha_i^*)^2\\
        =&-\frac{\theta}{p_i n}(\alpha_i^{(t)}-\alpha_i^*)^2+\frac{\theta}{p_in}((-\phi'_i(x_i^Tw^{(t)})-\alpha_i^*)^2 - (1- \frac{\theta}{p_i})(\kappa_i^{(t)})^2).
    \end{align}

    We can then obtain \eqref{eq: E[A]} by taking expectation on $i$ conditioned on $\alpha^{(t)}$. On the other side,
    \begin{align}
    &B^{(t+1)} - B^{(t)} \\
    =&\norm{w^{(t+1)}-w^*}^2 - \norm{w^{(t)}-w^*}^2\\
    =&\norm{w^{(t)} - \frac{\theta}{n \lambda p_i}\kappa_i^{(t)}x_i-w^*}^2 - \norm{w^{(t)}-w^*}^2\\
    =& \frac{2 \theta}{n \lambda p_i}(\kappa_i^{(t)}x_i)^T(w^{(t)}-w^*)+ \frac{\theta^2v_i}{n^2 \lambda^2 p_i^2}(\kappa_i^{(t)})^2.
    \end{align}

    From \eqref{eq: nabla(P)}, we obtain that $\E(\kappa_i^{(t)}/p_i) = n\nabla P(w^{(t-1)})$ and further take expectation on $i$ conditioned on $w^{(t)}$, \eqref{eq: E[B]} is verified.
\end{proof}



\begin{lemma}
    Assume that each $\phi_i$ is $L$-smooth and convex, then we have
    \begin{equation}\label{eq: strong ineq}
        \frac{1}{n}\sum_{i=1}^n\norm{\phi'_i(w^Tx_i) - \phi'_i((w^*)^Tx_i))}^2 \leq 2L(P(w)-P(w^*)- \frac{\lambda}{2}\norm{w-w^*}^2).
    \end{equation}
\end{lemma}

\begin{proof}
    We define
    \begin{equation}
        g_i(z) = \phi_i(z) - \phi_i(z^*) - \phi'_i(z^*)(z-z^*).
    \end{equation}

    Then we know $g_i$ is non-negative and $L$-smooth and further we have
    \begin{equation}
     \norm{g'_i(z)}^2 \leq 2 L g_i(z).
 \end{equation}

 Actually, note that $g_i$ is $L$-smooth, we have
 \begin{equation}
     g_i(z) - g_i(\hat{z}) \leq g'_i(\hat{z})(z-\hat{z}) + \frac{L}{2}(z-\hat{z})^2.
 \end{equation}

 Expected result is then given by setting $\hat{z} = z - g'_i(z)$.

 Therefore, we set $z = w^Tx_i$ and $z^* = (w^*)^Tx_i$, then
 \begin{equation}
     \norm{\phi'_i(w^Tx_i) - \phi'_i((w^*)^Tx_i))}^2 = \norm{g'_i(w^Tx_i)}^2 \leq 2Lg_i(w^Tx_i).
 \end{equation}

 Taking expectation on both side, we have
 \begin{align}
     &\frac{1}{n}\sum_{i=1}^n\norm{\phi'_i(w^Tx_i) - \phi'_i((w^*)^Tx_i))}^2 = \norm{g'_i(w^Tx_i)}^2 \\
     \leq ~&2L\E[\phi_i(w^Tx_i) - \phi_i((w^*)^Tx_i) - \phi'_i((w^*)^Tx_i)(w^Tx_i-(w^*)^Tx_i)]\\
     = ~& 2L(P(w)-\frac{\lambda}{2}\norm{w}^2-P(w^*) + \frac{\lambda}{2}\norm{w^*}^2 - \lambda (w^*)^T(w-w^*))\\
     = ~&2L(P(w)-P(w^*)- \frac{\lambda}{2}\norm{w-w^*}^2),
 \end{align}

 where we make use of
 \begin{equation}
    \E[\nabla P(w^*)] = \E[\phi'((w^*)^Tx_i) x_i + \lambda w^*] = 0.
\end{equation}
\end{proof}

Now we are able to prove Theorem \ref{thm: seperate convex} as follows.
\begin{proof}We set $\gamma = \lambda L$ and let $D^{(t)} = A^{(t)} + \gamma B^{(t)}$, then we have
\begin{align}
\E[D^{(t+1)}-D^{(t)}] &= \E[A^{(t+1)}-A^{(t)}] + \gamma\E[B^{(t+1)}-B^{(t)}]\\
&\overset{\eqref{eq: E[A]}, \eqref{eq: E[B]}}{=}-\theta A^{(t)} + \frac{\theta}{n}\sum_{i=1}^n((-\phi'_i(x_i^Tw^{(t)})-\alpha_i^*)^2 - (1- \frac{\theta}{p_i})(\kappa^{(t)}_i)^2) \\
&+ \gamma (-\frac{2 \theta }{\lambda}\nabla P(w^{(t)})^T(w^{(t)}-w^*) + \sum_{i=1}^n \frac{\theta^2v_i(\kappa^{(t)}_i)^2}{n^2 \lambda^2p_i})\\
&\overset{\eqref{eq: strong ineq}}{\leq} -\theta A^{(t)} + \underbrace{\sum_{i=1}^n(-\frac{\theta}{n}(1- \frac{\theta}{p_i})+\frac{\theta^2v_i \gamma}{n^2 \lambda^2 p_i})(\kappa_i^{(t)})^2}_{\mathbf{C}} \\
&+\theta(2L(P(w)-P(w^*)- \frac{\lambda}{2}\norm{w-w^*}^2)-\frac{2 \gamma}{\lambda}\nabla P(w^{(t)})^T(w^{(t)}-w^*))\\
&\overset{\gamma=\lambda L}{\leq} -\theta (A^{(t)} + \gamma\norm{w^{(t)}-w^*}^2) +  \mathbf{C}\\
&=-\theta D^{(t)} + \mathbf{C}.
\end{align}

The last second inequality comes from convexity of $P$, where we have $P(w^{(t)})-P(w^*) \leq \nabla P(w^{(t)})^T(w^{(t)}-w^*)$.
\end{proof}

Lemma \ref{lem: optimal probabilities} is easy to verify by derive KKT conditions of optimization problem \eqref{Prob: to derive best probabilities}.
\end{document}